\documentclass[a4paper,11pt]{amsart}
\usepackage{amsmath,amssymb}
\usepackage{enumerate}
\usepackage{graphics}

\title{Common hypercyclic vectors for composition operators}
\author{Fr\'ed\'eric Bayart}
\date{}

\address{
Universit\'e des Sciences et Technologies de Lille, Laboratoire Agat, U.F.R. de
Math\'ematiques Pures et Appliqu\'ees, B\^at M2., F-59655 Villeneuve
d'Ascq Cedex}
\email{bayart@agat.univ-lille1.fr}

\newcommand{\mtt}{\mathbb{T}}

\newcommand{\mtd}{\mathbb{D}}
\newcommand{\mtn}{\mathbb{N}}

\newcommand{\mtc}{\mathbb{C}}
\newcommand{\mtq}{\mathbb{Q}}
\newcommand{\mtr}{\mathbb{R}}

\newcommand{\cphi}{C_\phi}
\newcommand{\hdeux}{H^2(\mtd)}
\newcommand{\mch}{\mathcal{H}}

\newcommand{\mcb}{\mathcal{B}}

\newcommand{\mcl}{\mathcal{L}}
\newcommand{\dis}{\displaystyle}

\newcommand{\cplus}{\mtc_+}

\newcommand{\veps}{\varepsilon}
\newcommand{\vphi}{\varphi}

\newcommand{\cvphi}{C_\varphi}

\newcommand{\lpmu}{L^p\left(\mtr,\frac{dt}{(1+t^2)^\alpha}\right)}
\newcommand{\lpdmu}{L^p\left(\mtr,d\mu\right)}

\newcommand{\autd}{Aut(\mtd)}

\DeclareMathOperator{\supp}{supp}

\newtheorem{defi}{Definition}
\newtheorem{thm}{Theorem}
\newtheorem{prop}{Proposition}
\newtheorem{cor}{Corollary}
\newtheorem{lem}{Lemma}

\def\piednote#1{\let\oldfn=\thefootnote\def\thefootnote{}\footnote{\noindent#1}%
\addtocounter{footnote}{-1}\def\thefootnote{\oldfn}}

\newenvironment{pf}
  {\noindent\textbf{Proof :}}
  {\hfill $\blacksquare$\par\medskip}
 
\newenvironment{ex}
{ \noindent {\it Example\/} : }
{\null \hfill \par\medskip} 

\begin{document}


\begin{abstract}
We study the existence of a common hypercyclic vector for different families of composition operators. We also give a continuous version of Salas theorem on weighted shifts.
\end{abstract}

\maketitle

\setcounter{section}{-1}

\section{Introduction}
A continuous operator acting on a topological vector space $X$ is called \emph{hypercyclic} provided there exists a vector $x\in X$ such that its orbit $\{T^nx;\ n\geq 1\}$ is dense in $X$. Such a vector is called a hypercyclic vector for $T$. The set of hypercyclic vectors will be denoted by $HC(T)$. The first example of hypercyclic operator was given by Birkhoff, 1929 \cite{BIRKHOFFHYPER}, who shows that the operator of translation by a non-zero complex number is hypercyclic on the space of holomorphic functions. For a complete account on hypercyclicity, we refer to \cite{GROSSEERDMANN}.

The main focus of our study is the hypercyclic behavior for composition operators. Let us denote by $\hdeux$ the Hardy space on the unit disk $\mtd$, and by $\autd$ the set of automorphisms of $\mtd$. For $\vphi$ in $\autd$, the hypercyclicity of the composition operator $\cvphi$ defined on $\hdeux$ by $\cvphi(f)=f\circ\vphi$ is well-understood since the work of Bourdon and Shapiro \cite{BOURDONSHAPIMEM} :
\begin{thm}\label{THMBOURDONSHAPI}
$\cvphi$ is hypercyclic on $\hdeux$ if, and only if, $\vphi$ has no fixed point in $\mtd$.
\end{thm}
This theorem emphasizes a previous result of Seidel and Walsh \cite{SEIDELWALSH}, who proved the same theorem for $\cvphi$ acting on the space of holomorphic functions on $\mtd$.

We will concentrate on the common hypercyclicity of a family of operators. Given a family $(T_\lambda)_{\lambda\in\Lambda}$ of hypercyclic operators on $X$, we ask whether it is possible to find a single vector $x$ which is hypercyclic for all $T_\lambda$. Observe that if the family is countable, and if $X$ is a F-space, a Baire's categorical argument implies that this is always possible : indeed, it turns out that $HC(T)$ is either empty or a dense $G_\delta$ set. For uncountable families, the first positive answer was given by E.Abakumov and J.Gordon \cite{ABAKGORD}, emphazing a theorem of Rolewicz :
\begin{thm}\label{THMABAKGORDBACK}
Let $B$ be the backward shift acting on $\ell^2$. There exists a common hypercyclic vector for the operators $\lambda B$, $\lambda>1$. 
\end{thm}

In section 1, we will recall as a theorem the construction made in the paper of Abakumov and Gordon. We will deduce a criterion for common hypercyclicity of multiples of a single operator, and we will apply this criterion to adjoints of multipliers. Section 2 is devoted to some positive and negative results for the problem of simultaneous hypercyclicity of composition operators. In particular, theorem \ref{THMSEIDELWALSHSIMUL} below is a simultaneous version of the theorem of Seidel/Walsh. Let us mention that the situation here is more complicated than in Birkhoff's theorem, since you have to handle not only translations, but also homotheties. Finally, in section 3, we provide some remarks and problems. In particular, we give a continuous analog to some well-known theorems on weighted shifts.\\
\noindent\textbf{Acknowledgements.} We thank E.Abakumov and J.Saint-Raymond for their help.

\section{Adjoints of multipliers}

\subsection{The size of the set of common hypercyclic vectors}

We begin by the following result, suggested by J.Saint-Raymond (the same was used in \cite[sec. 3.4]{ABAKGORD} :
\begin{prop}
Let $X$ be a $F-$space, $A\subset \mcl(X)$ such that $A$ is the countable union of compact sets. Then $\dis\bigcap_{T\in A}HC(T)$ is a $G_\delta$ set.
\end{prop}
\begin{pf}
Define $M=\big\{(T,x)\in A\times X;\ x\notin HC(T)\big\}$, and consider $\big(\mcb_m\big)$ a countable basis of open sets in $X$. Then :
\begin{eqnarray*}
M^c&=&\big\{(T,x)\in A\times X;\ x\in HC(T)\big\}\\
&=&\bigcap_{m\geq 1}\bigcup_{n\geq 0}\big\{(T,x);\ T^nx\in\mcb_m\big\}.
\end{eqnarray*}
In particular, $M^c$ is a $G_\delta$ set in $A\times X$. Let us write $M=\bigcup_{k\geq 1}F_k$ (resp. $A=\bigcup_{p\geq 1}A_p$) where each $F_k$ is closed in $A\times X$ (resp. each $A_p$ is compact). If $\pi$ denotes the projection of $\mcl(X)\times X\to X$ onto the second coordinate, we deduce that :
$$\pi(M)=\pi\left(\bigcup_{k\geq 1}F_k\right)=\bigcup_{k\geq 1}\bigcup_{p\geq 1}\pi\big(F_k\cap(A_p\times X)\big).$$
Each set $\pi\big(F_k\cap(A_p\times X)\big)$ is closed in $X$ since $A_p$ is compact and $F_k$ is closed. Therefore, $\pi(M)$ is $F_\sigma$. Now, $\pi(M)=\left[\dis\bigcap_{T\in A}HC(T)\right]^c$, and this gives the conclusion.
\end{pf}
The previous proposition does not ensure that $\dis\bigcap_{T\in A}HC(T)$ is not empty. But as soon as this is the case, we should control the size of this set :
\begin{cor}\label{CORNATUREHYPERCYC}
Let $X$ be a F-space, $A\subset\mcl(X)$. Assume that :
\begin{enumerate}
\item $A$ is the countable union of compact sets.
\item $\dis\bigcap_{T\in A}HC(T)\neq\varnothing$.
\item There exists $S\in A$ which commutes with all $T\in A$.
\end{enumerate}
Then $\dis\bigcap_{T\in A}HC(T)$ is residual (dense $G_\delta$).
\end{cor}
\begin{pf}
Pick $x\in\dis\bigcap_{T\in A}HC(T)$, and $S$ as in 3. It is straightforward that the dense set $\left\{S^kx;\ k\geq 1\right\}$ is contained in $\dis\bigcap_{T\in A}HC(T)$.
\end{pf}
\subsection{Abakumov-Gordon's construction}
Our proofs will be constructive ones. We will need an approximation tool, provided by the paper of Abakumov and Gordon :
\begin{thm}\label{THMABAKGORDMULT}
There exists an integer $k_0\geq 1$, a function $j:\{n\in\mtn;\ n\geq k_0\}\to\mtn$ such that, for any sequence $(\alpha_l)_{l\geq 1}$ of positive real numbers, there exists a sequence $(M_k)_{k\geq k_0}$ of positive integers, a sequence $(r_k)_{k\geq k_0}$ of positive real numbers such that :
\begin{enumerate}[i.]
\item $(M_k)$ is increasing, $M_{k+1}-M_k\to+\infty$.
\item $(r_k)$ is decreasing, $\frac{r_{k+1}}{r_k}\to 0$.
\item For any $l$ of $\mtn$, for any $\varepsilon>0$, for any $\lambda>1$, for any $K>0$, there exists $k>K$ such that :
$$j(k)=l\textrm{ and }|\lambda^{M_k}r_k-\alpha_l|<\varepsilon.$$
\end{enumerate}
\end{thm}
$j$ is a \emph{choice function}. This theorem can be seen as an uncountable Baire's type theorem. It is trivial that :
$$\forall \lambda>1,\ \exists (M_k)_{k\in\mtn}\in\mtn,\ (r_k)_{k\in\mtn}\in\mtr\textrm{ such that }\{\lambda^{M_k}r_k\}\ \textrm{is dense in $\mtr_+$}.$$
Theorem \ref{THMABAKGORDMULT} says :
$$\exists(M_k)_{k\in\mtn}\in\mtn,\ (r_k)_{k\in\mtn}\in\mtr\textrm{ such that }\forall \lambda>1,\ \{\lambda^{M_k}r_k\}\ \textrm{is dense in $\mtr_+$}.$$
We will also need an additive version of this theorem :
\begin{thm}\label{THMABAKGORDADD}
There exists an integer $k_0\geq 1$, a function $j:\{n\in\mtn;\ n\geq k_0\}\to\mtn$,  a sequence $(M_k)_{k\geq k_0}$ of positive integers, a sequence $(X_k)_{k\geq k_0}$  of real numbers such that :
\begin{enumerate}
\item $(M_k)$ is increasing, $M_{k+1}-M_k\to+\infty$.
\item $(X_k)$ is increasing,  $X_{k+1}-X_k\to+\infty$.
\item For any $l$ in $\mtn$, for any $\veps>0$, for any $a>0$, for any $K>0$, there exists $k>K$ such that :
$$j(k)=l,\textrm{ and }\big|M_ka-X_k\big|<\veps.$$
\end{enumerate}
\end{thm}
\subsection{A criterion for common hypercyclicity}
When one wants to show that an operator $T$ is hypercyclic, the most useful tool is the hypercyclic criterion formulated first by C.Kitai (see \cite[cor. 1.5]{GODSHAPI} for a statement of this criterion). We give here a sufficient condition for the existence of a common hypercyclic vector for all multiples of an operator.
\begin{thm}\label{THMKITAISIMUL}
Let $X$ be a separable Banach space, and $T\in L(X)$. Assume that :
\begin{enumerate}[a)]
\item $V=\bigcup_n Ker(T^n)$ is dense in $X$.
\item There exists $S:V\to X$ with $TS=Id_V$ and $\|Sx\|\leq \|x\|$ for all $x$ in $V$.
\end{enumerate}
Then $\dis\bigcap_{\lambda>1}HC(\lambda T)$ is dense $G_\delta$.
\end{thm} 
\begin{pf}
By corollary \ref{CORNATUREHYPERCYC}, it is enough to prove that $\bigcap_{\lambda>1}HC(\lambda T)$ is non-empty. Fix $(v_l)$ a dense sequence in $V$, and set $\alpha_l=\|v_l\|$. We define the function $j$ and the sequences $(M_k)$ and $(r_k)$ as in theorem \ref{THMABAKGORDMULT}. For $k\geq k_0$, let us set :
\begin{itemize}
\item $d_k=r_k-r_{k+1}\geq 0$. 
\item $w_k=v_{j(k)}$ if $T^{M_{k+1}-M_k}v_{j(k)}=0$, $w_k=0$ otherwise.
\item $\dis y_k=\frac{d_k}{\|w_k\|}S^{M_k}w_k$ if $w_k\neq 0$, $y_k=0$ otherwise.
\end{itemize}
We claim that $f=\sum_{m\geq k_0}y_m$ is hypercyclic for each $\lambda T$, with $\lambda>1$. First, observe that if $m<k$, $T^{M_k}S^{M_m}w_m=T^{M_k-M_m}w_m=0$, what implies :
$$\left\|T^{M_k}f\right\|=\left\|\sum_{m\geq k}\frac{d_m}{\|w_m\|}T^{M_k}S^{M_m}w_m\right\|\leq\sum_{m\geq k}d_k=r_k.$$
Take now $\veps>0$ and $l\in\mtn$. By theorem \ref{THMABAKGORDMULT}, there exists $k\in\mtn$ such that :
\begin{itemize}
\item $j(k)=l$, and $w_k=v_l$.
\item $\left|\lambda^{M_k}r_k-\|v_l\|\right|\leq\varepsilon$.
\item $\frac{r_{k+1}}{r_k}(\varepsilon+\|v_l\|)\leq\varepsilon$.
\end{itemize}
Then,
\begin{eqnarray*}
\left\|(\lambda T)^{M_k}f-v_l\right\|&\leq&\left\|\lambda^{M_k}\frac{d_k}{\|v_l\|}v_l-v_l\right\|+\left\|\sum_{m>k}\lambda^{M_k}\frac{d_m}{\|w_m\|}T^{M_k}S^{M_m}w_m\right\|\\
&\leq&\|v_l\|\left(\left|\frac{\lambda^{M_k}r_k}{\|v_l\|}-1\right|+\frac{\lambda^{M_k}r_{k+1}}{\|v_l\|}\right)+\lambda^{M_k}r_{k+1}\\
&\leq&\varepsilon+2\lambda^{M_k}r_{k+1}\\
&\leq&\varepsilon+2(\varepsilon+\|v_l\|)\frac{r_{k+1}}{r_k}\\
&\leq&3\varepsilon.
\end{eqnarray*}
This achieves to prove that $f$ is hypercyclic for $\lambda T$.
\end{pf}
Hypercyclic operators are strongly connected with the existence of invariant subspaces. The following corollary illustrates this link :
\begin{cor}
Under the assumptions of theorem \ref{THMKITAISIMUL}, there exists a dense subspace of $X$, invariant by $T$, whose elements, except 0, are hypercyclic vectors for $\lambda T$, with $\lambda>1$.
\end{cor}
\begin{pf}
Take $x$ in $\bigcap_{\lambda>1}HC(\lambda T)$. Then
$$M=\{p(T)x;\ p\textrm{ is a polynomial}\}$$
answers the question : the proof given by P.S. Bourdon in \cite{BOURDONPAMS} works also in this setting.
\end{pf}
\subsection{Application to adjoints of multipliers}
\begin{cor}\label{CORSIMULADJOINTS}
Let $\vphi$ be an inner function, not a constant, and $M_\vphi$ be the associated multiplier on $\hdeux$ (defined by $M_\vphi(f)=\vphi f$). Then $\dis\bigcap_{\lambda>1}HC(\lambda M_\vphi^*)$ is a residual set. 
\end{cor}
By choosing $\vphi(z)=z$, we retrieve theorem \ref{THMABAKGORDBACK}.\\
\begin{pf}
\begin{enumerate}[a)]
\item It is plain that $\ker \left(M_\vphi^*\right)^n=\left(\vphi^n H^2\right)^\perp$. Let us recall the following result from \cite[p. 34-35]{NIKOLSKIBOOK} : let $E$ be a normed space, and $(E_n)$ be a sequence of subspaces of $E$. We define :
$$\underline{\lim}\ E_n=\big\{x\in E;\ \lim_n dist(x,E_n)=0\big\}.$$
If $E=H^2$, and if $E_n=\left(\theta_n H^2\right)^\perp$, where $(\theta_n)$ is a sequence of inner functions, then we have :
\begin{lem}\label{LEMNIKOLSKI}
$$\underline{\lim}\ \big(\theta_n\hdeux)^\perp=H^2(\mtd)\iff \forall z\in\mtd,\ \lim_n\theta_n(z)=0.$$
\end{lem}
In our context, $\theta_n=\vphi^n$, and $\left(\vphi^n H^2\right)^\perp\subset \left(\vphi^{n+1} H^2\right)^\perp$. So, $\underline{\lim}(\vphi^nH^2)^\perp\subset\bigcup_n(\vphi^nH^2)^\perp$. Now, since $\vphi$ is not constant, for each $z$ in $\mtd$, $\vphi^n(z)\to 0$, and lemma \ref{LEMNIKOLSKI} gives 
$$\hdeux\subset\overline{\bigcup_n Ker\big((M_\vphi^*)^n\big)}.$$
\item If $f\in V$, and $g\in \hdeux$, then :
\begin{eqnarray*}
<g,M_{\varphi}^*M_\varphi f>&=&<M_\varphi g,M_\varphi f>\\
&=&<g,f> \textrm{ since $\varphi$ is inner.}
\end{eqnarray*}
So we can take $S=M_\vphi$ in part b) of theorem \ref{THMKITAISIMUL}.
\end{enumerate}
\end{pf}
\section{Composition operators}
\subsection{Geometry of the disk}
For details on the background material of this section, we refer to \cite{SHAPIROBOOK}. The automorphisms of $\mtd$ can be classified in function of their fixed points : $\vphi\in \autd$ is called :
\begin{itemize}
\item \emph{parabolic} if $\vphi$ has a single (attractive) fixed point on $\mtt=\partial\mtd$.
\item \emph{hyperbolic} if $\vphi$ has an attractive fixed point on $\mtt$, and a second one on $\mtt$.
\item \emph{elliptic} if $\vphi$ has an attractive fixed point in $\mtd$.
\end{itemize}
We are concerned by parabolic and hyperbolic automorphisms. It is easier to describe their action on the right half-plane $\cplus$. Denote $\sigma:\mtd\to\cplus$, $\sigma(z)=\frac{1+z}{1-z}$ the Cayley map from $\mtd$ onto $\cplus$. For $\vphi\in\autd$ with $+1$ as attractive fixed point, set $\psi=\sigma\circ\vphi\circ\sigma^{-1}$. Then :
\begin{itemize}
\item $\psi(z)=z+ia$ where $a\in\mtr$, $a\neq 0$, if $\vphi$ is parabolic (a parabolic automorphism of $\mtd$ is conjugated to a translation).
\item $\psi(z)=\lambda(z-ib)+ib$, where $\lambda>1$ and $b\in\mtr$, if $\vphi$ is hyperbolic (a hyperbolic automorphism of $\mtd$ is conjugated to a positive dilatation).
\end{itemize}
\subsection{Main statements} In view of theorem \ref{THMBOURDONSHAPI}, it comes a natural \mbox{question :}
\begin{quote}
Does there exist a common hypercyclic vector for all composition operators $\cvphi$ on $\hdeux$, where $\vphi\in\autd$ has no fixed point in $\mtd$?
\end{quote}
Here, you can play with two parameters : you can choose the attractive fixed point, and its attractivity (the scalars $\lambda$, $a$ and $b$ of the previous paragraph). The following result shows that it is impossible to have a wide set of attractive fixed points :
\begin{thm}
Let $A\subset\autd$ such that, for any $\vphi$ in $A$, $\vphi$ has no fixed point in $\mtd$. Moreover, suppose that :
$$B=\left\{\omega\in\mtt;\ \exists\vphi\in A\textrm{ such that $\omega$ is the attractive fixed point of $\vphi$}\right\}$$
has positive Lebesgue measure. Then $\bigcap_{\vphi\in A}HC(\cvphi)=\varnothing$.
\end{thm}
Here, $\cvphi$ is considered as a composition operator on $\hdeux$.\\
\begin{pf}
The theorem is a direct consequence of the following lemma, since a function of $\hdeux$ admits angular limits almost everywhere on the boundary.
\begin{lem}
Suppose that $\vphi\in\autd$, and that $\omega\in\mtt$ is the attractive fixed point of $\vphi$. If $f\in\hdeux$ is a hypercyclic vector for $\cvphi$, then $f$ has no angular limit at $\omega$.
\end{lem}
\begin{pf}
We denote $\vphi_n=\vphi\circ\dots\circ\vphi$ ($n$ times). By Denjoy-Wolff's theorem, $(\vphi_n(0))$ converges nontangentially to $\omega$. Now, evaluation at 0 is continuous on $\hdeux$, and by hypercyclicity of $f$, there exist integers $m$ and $n$, as large as possible, such that :
$$\left|f\circ\vphi_m(0)-0\right|<1/4\textrm{ and }\left|f\circ\vphi_n(0)-1\right|<1/4.$$
In particular, $f$ does not admit any nontangential limit at $\omega$.
\end{pf}
\end{pf}
So, essentially we have to fix the attractive fixed point, say +1, and the question becomes :
\begin{quote}
Does there exist a common hypercyclic vector for all composition operators $\cvphi$ on $\hdeux$, where $\vphi\in\autd$ has +1 as attractive fixed point?
\end{quote}
We are not able to give a positive or a negative answer to this question. But if we relax the hypothesis, this will be the case. On the one hand, we can forget the growth condition : if $\vphi\in\autd$, $\cvphi$ is a composition operator on $H(\mtd)$, the $F-$space of holomorphic functions on $\mtd$. Since the theorem of Seidel/Walsh, we know that such a composition operator is hypercyclic. Under these assumptions, there exists a common hypercyclic vector :
\begin{thm}\label{THMSEIDELWALSHSIMUL}
Let $\omega\in\mtt$. There exists a common hypercyclic vector for all composition operators $\cvphi$ acting on $H(\mtd)$, where $\vphi\in\autd$ admits $\omega$ as attractive fixed point. Moreover, the set of common hypercyclic vectors is a residual set.
\end{thm}
On the other hand, we can ignore the regularity condition : by results of Nordgren \cite{NORDGREN}, $\cvphi$ is also a composition operator on $L^2(\mtt)$. An application of Kitai's criterion should prove its hypercyclicity. We directly prove a simultaneous hypercyclicity theorem :
\begin{thm}\label{THMSIMULL2}
Let $\omega\in\mtt$. There exists a common hypercyclic vector for all composition operators $\cvphi$ acting on $L^2(\mtt)$, where $\vphi\in\autd$ admits $\omega$ as attractive fixed point. Moreover, the set of common hypercyclic vectors is a residual set.
\end{thm}
The remaining part of this section is devoted to the proof of the previous theorems. We will assume that $\omega=+1$.
\subsection{Proof of the holomorphic case}
We take the model of the half-plane. Define $T_a$ and $S_{\lambda,b}$ by :
$$T_a(f)(z)=f(z+ia)$$
$$S_{\lambda,b}(f)(z)=f(\lambda(z-ib)+ib).$$
It suffices to show that $\bigcap_{a\neq 0}HC(T_a)$ and $\bigcap_{\lambda>1,b\in\mtr}HC(S_{\lambda,b})$ are dense $G_\delta$. Till the end of this section, we fix $(\delta_k)$, $0<\delta_k<1$ a sequence which converges to 0, and $(P_l)$ a sequence in $H(\cplus)$ such that, for any $\mu\geq 1$ and any $\tau\in\mtr$, $\big(P_l(\mu z-\mu i\tau)\big)$ is dense in $H(\cplus)$ (for example, $(P_l)$ could be the sequence of polynomials with coefficients in $\mtq$). We handle separately the parabolic and the hyperbolic case.
\subsubsection{Parabolic automorphisms}
By corollary \ref{CORNATUREHYPERCYC}, it is enough to prove by instance that $\bigcap_{a>0}HC(T_a)$ is not empty. We fix sequences $(M_k)$ and $(X_k)$ as in theorem \ref{THMABAKGORDADD}. For $k\geq k_0+1$, let us set :
$$R_k=\min\left(\frac{X_{k+1}-X_k}{2},\frac{X_k-X_{k-1}}{2}\right).$$
We build by induction rectangles $C_k$, $D_k$ and $\Gamma_k$, for $k\geq k_0+1$, beginning by the initialization $\Gamma_{k_0}=\left\{(1,0)\right\}$. For $k\geq k_0+1$, fix $C_k$ the square whose center is $(R_k/2,0)$, whose side has length $R_k-\delta_k$. Observe that, for any compact $K$ of $\cplus$, for $k$ large enough, $K\subset C_k$. Consider $D_k=C_k+iX_k$. The squares $(D_k)$ are disjoint. Moreover, there exists a rectangle $\Gamma_k$ which contains $\Gamma_{k-1}$, $D_k$, but which has empty intersection with $D_{k+1}$.

We then define a sequence $(\pi_k)_{k\geq k_0}$ of polynomials, by setting $\pi_{k_0}(z)=1$. Next, for $k>k_0$, Runge's theorem gives a polynomial $\pi_k$ \mbox{satisfying :}
\begin{enumerate}[a)]
\item $\left|\pi_k(z)-P_l(z-iX_k)\right|\leq\frac{1}{2^k}$ if $z\in D_k$ and $j(k)=l$.
\item $\left|\pi_k(z)-\pi_{k-1}(z)\right|\leq \frac{1}{2^k}$ if $z\in\Gamma_{k-1}$.
\end{enumerate}
The sequence $(\pi_k)$ converges uniformly on each compact subset of $\cplus$. Let us denote by $f$ its limit. Observe that, for each $z\in\Gamma_k$, we have :
$$|f(z)-\pi_k(z)|\leq |\pi_k(z)-\pi_{k+1}(z)|+|\pi_{k+1}(z)-\pi_{k+2}(z)|+\dots\leq\frac{1}{2^k}.$$
We claim that $f$ is hypercyclic for each $T_a$, with $a>0$. Indeed, fix $l\in\mtn$, $K$ a compact subset of $\cplus$, and $\eta>0$ such that $K_1=K+\overline{B}(0,\eta)\subset\cplus$. Let $0<\delta<\eta$ with :
$$z_1,z_2\in K_1\ \wedge\ |z_1-z_2|\leq\delta\implies|P_l(z_1)-P_l(z_2)|\leq\varepsilon.$$
There exists an integer $k$ such that :
\begin{itemize}
\item $j(k)=l$, $\frac{1}{2^k}\leq\veps$, and $K_1\subset C_k$.
\item $|aM_k-X_k|\leq\delta$.
\end{itemize}
Then, for $z\in K$, $z+iM_ka-iX_k\in K_1\subset C_k$, and therefore $z+iM_ka\in D_k$. This implies that :
\begin{eqnarray*}
\left|[T_a(f)]^{M_k}(z)-P_l(z)\right|&\leq&\veps+\left|\pi_k(z+iM_ka)-P_l(z)\right|\\
&\leq&2\veps+\left|P_l(z+iM_ka-iX_k)-P_l(z)\right|\\
&\leq&3\veps.
\end{eqnarray*}
\subsubsection{Hyperbolic automorphisms}
\begin{figure}[!t]
\begin{center}
\input{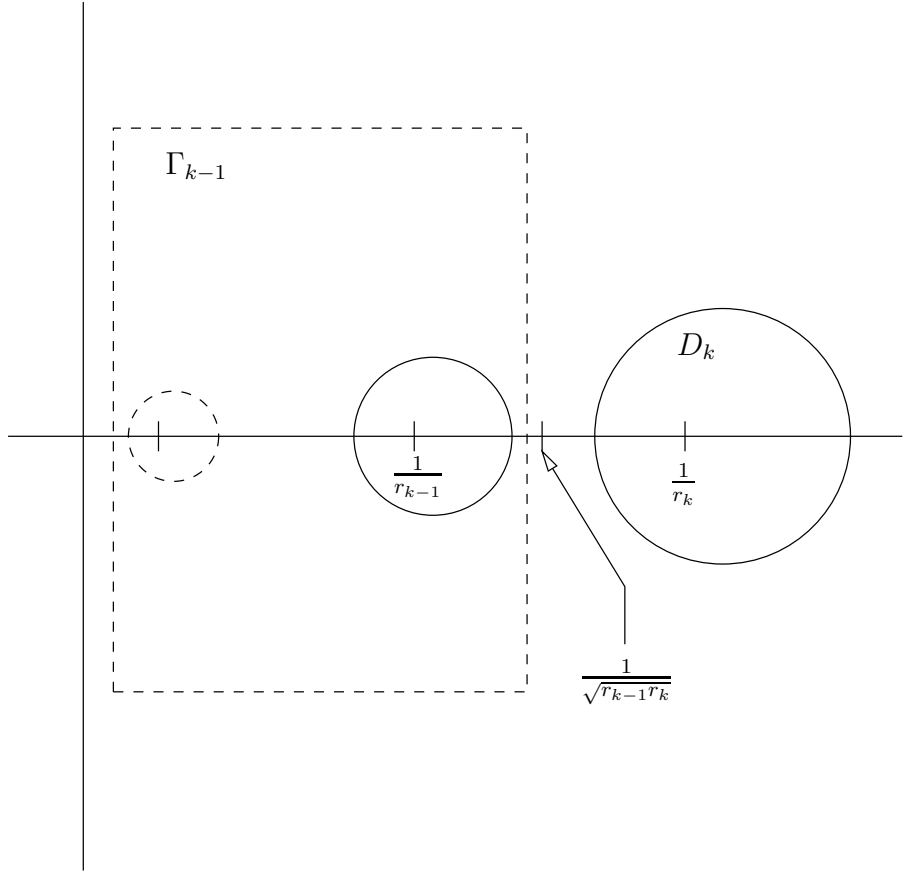}
\end{center}
\caption{The hyperbolic construction\label{FIGHYPER}}
\end{figure}

Here, dilatations do not commute, and we need to prove that $\bigcap_{\lambda,b}HC(S_{\lambda,b})$ is dense. First, apply theorem \ref{THMABAKGORDMULT} with the sequence $(\alpha_l)$ identically one, to obtain sequences $(M_k)$ and $(r_k)$. For $k\geq k_0+1$, we set :
$$R_k=\delta_k\min\left(\frac{\sqrt{\frac{r_{k-1}}{r_k}}-1}{\sqrt{\frac{r_{k-1}}{r_k}}+1},\frac{1-\sqrt{\frac{r_{k+1}}{r_k}}}{\sqrt{\frac{r_{k+1}}{r_k}}+1}\right).$$
We always fix $\Gamma_{k_0}=\left\{(1,0)\right\}$, and for $k\geq k_0$, we consider $C_k$ the hyperbolic disk whose center is $(1,0)$ and whose radius is $R_k$ :
$$C_k=\left\{z\in\mtc_+;\ \frac{|z-1|}{|z+1|}\leq R_k\right\}.$$
Let $D_k$ be the image of $C_k$ by the homothety of center 0, of ratio $\frac{1}{r_k}$. $(D_k)$ are disjoint sets, and by construction there exists a rectangle $\Gamma_k$ which contains $\Gamma_{k-1}$ and $D_k$, but whose intersection with $D_{k+1}$ is empty (see figure \ref{FIGHYPER}).
Finally, we set $\pi_{k_0}(z)=1$, and if $k>k_0$, $l=j(k)$, Runge's theorem gives us   a polynomial $\pi_k$ which satisfies :
\begin{enumerate}[a)]
\item $|\pi_k(z)-P_l(r_k z)|\leq\frac{1}{2^k}$ if $z\in D_k$.
\item $|\pi_k(z)-\pi_{k-1}(z)|\leq\frac{1}{2^k}$ if $z\in\Gamma_{k-1}$.
\end{enumerate}
As previously, $(\pi_k)$ converges uniformly on each compact of $\cplus$ to a function $f$, with 
$$|f(z)-\pi_k(z)|\leq\frac{1}{2^k}\textrm{ if }z\in\Gamma_k.$$
For $\mu>1$, we claim that $g(z)=f(\mu z)$ is hypercyclic for each $S_{\lambda,b}$, $\lambda>1$, $b\in\mtr$. Indeed, fix $l\in\mtn$, $\veps>0$, $K$ a compact of $\cplus$ and $\eta>0$ such that $K_1=K+\overline{B}(0,\eta)\subset\cplus$. Let $0<\delta<\eta$ with :
$$z_1,z_2\in K_1\ \wedge\ |z_1-z_2|\leq\delta\implies|P_l(\mu z_1-\mu ib)-P_l(\mu z_2-\mu ib)|\leq\varepsilon.$$
By theorem \ref{THMABAKGORDMULT}, there exists an integer $k$ with :
\begin{itemize}
\item $j(k)=l$ and $\frac{1}{2^k}\leq\veps$.
\item $\mu\lambda^{M_k}r_k(K-ib)+\mu r_kib\subset C_k$.
\item By letting $M$ such that $z\in K\implies |z|\leq M$, then :
$$\mu|\lambda^{M_k}r_k-1|(M+|b|)+\mu r_k|b|<\delta.$$
\end{itemize}
Then, if $z\in K$, $\mu\lambda^{M_k}(z-ib)+\mu ib\in D_k\subset\Gamma_k$, and so : 
\begin{eqnarray*}
|\left[S_{\lambda,b}(g)\right]^{M_k}(z)-P_l(\mu z-\mu ib)|
&=&\left|f\big(\mu\lambda^{M_k}(z-ib)+\mu ib\big)-P_l(\mu z-\mu ib)\right|\\
&\leq&\varepsilon+|\pi_k(\mu\lambda^{M_k}(z-ib)+\mu ib)-P_l(\mu z-\mu ib)|\\
&\leq&2\varepsilon+|P_l(\mu \lambda^{M_k}r_k(z-ib)+\mu r_kib)-P_l(\mu z-\mu ib)|\\
&\leq&3\varepsilon,
\end{eqnarray*}
where the last inequality comes from :
\begin{eqnarray*}
\left|\mu\lambda^{M_k}r_k(z-ib)+\mu r_kib-\mu z-\mu ib\right|&\leq&\mu\left|\lambda^{M_k}r_k-1\right|\big(|z|+|b|\big)+\mu r_k|b|\\
&<&\delta.
\end{eqnarray*}
Therefore, $\left\{f(\mu z);\ \mu\geq 1\right\}\subset\bigcap_{\lambda>1,b\in\mtr}HC(S_{\lambda,b})$, and $\left\{f(\mu z);\ \mu\geq 1\right\}$ is dense in $H(\cplus)$ since for example $f$ is hypercyclic for $S_{2,0}$.
\subsection{Proof of the $L^2-$case}
Let $\lambda_i$ be the probability measure on $\mtr$ defined by $d\lambda_i(t)=\pi^{-1}(1+t^2)^{-1}dt$ ($\lambda_i$ is the image of the Lebesgue measure on $\mtt$ by $\sigma$). Notice that $f\in L^2(\mtt)\iff f\circ\sigma^{-1}\in L^2(\mtr,d\lambda_i)$, and that :
$$\int_{-\infty}^{+\infty}|f\circ\sigma^{-1}(it)|^2d\lambda_i(t)=\frac{1}{2\pi}\int_{-\pi}^{\pi}|f(e^{i\theta})|^2d\theta.$$ 
Let us change the notations to avoid the integration on $i\mtr$. For $\lambda>1$, $a\in \mtr$  $(a\neq 0)$ and $b\in\mtr$, we now set :
$$T_a(f)(x)=f(x+a)$$
$$S_{\lambda,b}(f)(x)=f\big(\lambda(x-b)+b\big).$$
We prove the slightly more precise result :
\begin{thm}\label{THMSIMULPOIDS}
Let $p\geq 1$, $\alpha>1/2$, and consider $T_a$ and $S_{\lambda,b}$ as operators on $\lpmu$. Then $\bigcap_{a\neq 0}HC(T_a)$ and $\bigcap_{\lambda>1,b\in\mtr} HC(S_{\lambda,b})$ are dense $G_\delta$ in $\lpmu$.
\end{thm}
The case $p=2$ and $\alpha=1$ corresponds exactly to theorem \ref{THMSIMULL2}. The following lemma will be useful for our purpose :
\begin{lem}\label{LEMSUITEUTILE}
Let $(v_k)_{k\geq 1}$ be a non-decreasing sequence of positive numbers, which tends to $+\infty$. Then there exists a non-decreasing sequence $(u_k)_{k\geq 1}$ of positive numbers, which tends to $+\infty$, and such that :
\begin{enumerate}[a)]
\item $\dis \sum_{k\geq 1}\frac{u_k}{k^3}<+\infty$.
\item $\dis \frac{u_k}{v_k^3}\xrightarrow{k\to+\infty} 0$.
\item $\dis \sum_{m>k}\frac{u_m}{\left( (m-k)+v_m\right)^3}\xrightarrow{k\to+\infty} 0$.
\end{enumerate}
\end{lem}
\begin{pf}
For $k\geq 1$, we set $u'_k=\inf (k,v_{[k/2]},v_{[k/2]+1},\dots,v_k)$, and $u_k=\inf_{l\geq k} u_l'$. Assertions $a)$ and $b)$ are trivial. For $c)$ :
\begin{itemize}
\item $\dis \sum_{m>2k}\frac{u_m}{\left( (m-k)+v_m\right)^3}=\sum_{m>k}\frac{u_{m+k}}{(m+v_{m+k})^3}\leq \sum_{m>k}\frac{1}{m^2}\to 0$.
\item $\dis \sum_{k<m\leq 2k}\frac{u_m}{\left((m-k)+v_m\right)^3}\leq v_k\sum_{m\leq k}\frac{1}{(m+v_k)^3}\leq \frac C{v_k}\to 0$.
\end{itemize}
\end{pf}
\subsubsection{Parabolic automorphisms}
First, we prove that $\bigcap_{a>0}HC(T_a)$ is not empty (and therefore is a dense $G_\delta$) for $\alpha=2$. We set $d\mu=\frac{dt}{(1+t^2)^2}$, and $C>0$ is a constant such that, for $x>0$, $\dis \int_x^{+\infty} d\mu\leq \frac{C}{x^3}$. We consider sequences $(M_k)$, $(X_k)$ as in theorem \ref{THMABAKGORDADD}. In particular, we will assume that $X_k\geq k$. For $k\geq k_0$, let us define $\dis R_k=\inf\left(\frac{X_{k+1}-X_k}{2},\frac{X_k-X_{k-1}}{2},\frac{X_k}{2}\right)$. Without lost of generality, we can always assume that $(R_k)$ is increasing. Next, $(u_k)$ is defined by applying lemma \ref{LEMSUITEUTILE} to the sequence $(v_k)$ with $v_k=R_k-2$. $(f_l)$ is a dense sequence in $\lpdmu$ of compactly supported bounded functions, with $\|f_l\|_\infty^p\leq u_l$. For $k\geq k_0$ and $l=j(k)$, let us set :
\begin{itemize}
\item $w_k=f_l$ if $\supp\ f_l\subset[-R_k;R_k]$, and $w_k=0$ otherwise.
\item $h_k(x)=w_k(x-X_k)$ : $(h_k)$ have mutually disjoint supports.
\end{itemize}
Finally we define $f=\sum_{k\geq k_0}h_k$. First of all, $f\in\lpdmu$. Indeed,
$$\|f\|_p^p\leq \sum_{k\geq k_0}\int_{X_k/2}^{+\infty}|w_k(x-X_k)|^p d\mu\leq C\sum_{k\geq k_0}\frac{2^3u_k}{X_k^3}<+\infty.$$
We claim that $f$ is hypercyclic for $T_a$, with $a>0$. Indeed, let $l\in\mtn$, $\veps>0$ and $0<\delta<1$ whose value will be precised later. There exists $k\geq k_0$, as large as necessary, such that :
\begin{itemize}
\item $j(k)=l$ and $\supp\ f_l\subset[-R_k,R_k]$.
\item $|M_ka-X_k|\leq \delta$.
\item For $m\geq k$, $X_{m+1}-X_m\geq 1$.
\end{itemize}
Then,
$$
\big\|T_a^{M_k}f-f_l\big\|_p\leq \big\|T_a^{M_k}h_k-f_l\big\|_p+\big\|\sum_{m>k} T_a^{M_k}h_m\big\|_p+\big\|\sum_{m<k}T_a^{M_k}h_m\big\|_p.
$$
Now,
\begin{enumerate}[1.]
\item $\dis\|T_a^{M_k}h_k-f_l\|_p=\|T_{M_ka-X_k}f_l-f_l\|_p\leq\varepsilon$ as soon as $\delta$ is small enough.
\item \begin{eqnarray*}
 \left\|\sum_{m>k}T_a^{M_k}h_m\right\|_p^p&\leq& \sum_{m>k}\int_{X_m-R_m-M_k a}^{+\infty}|h_m(x+M_ka)|^pd\mu\\
&\leq& C\sum_{m>k}\frac{u_m}{(X_m-X_k-R_m-1)^3}.
\end{eqnarray*} Observe that :
$$X_m-X_{k}\geq X_m-X_{m-1}+\dots+X_{k+1}-X_k\geq 2R_m+m-k-1.$$
We deduce that :
$$\left\|\sum_{m>k}T_a^{M_k}h_m\right\|_p^p\leq C\sum_{m>k}\frac{u_m}{\left((m-k)+R_m-2\right)^3},$$
and this last quantity is smaller than $\varepsilon$ if $k$ is large enough.
\item $$
\begin{array}{rcll}
\dis \left\|\sum_{m<k}T_{M_ka}h_m\right\|_p^p&\leq&\dis \sum_{m<k}u_k\int_{-M_ka+X_m-R_m}^{-M_ka+X_m+R_m}d\mu\\
&\leq&\dis C\frac{u_k}{(M_ka-X_{k-1}-R_{k-1})^3}&\textrm{ (disjoint supports)}\\
&\leq&\dis C\frac{u_k}{(X_k-X_{k-1}-R_{k-1}-1)^3}\\
&\leq&\dis C\frac{u_k}{(R_k-1)^3}\xrightarrow{k\to+\infty} 0.
\end{array}$$
\end{enumerate}
It remains to prove the case $\alpha\neq 2$. We use a classical lemma (see \cite[p111 ``The hypercyclic comparison principle'']{SHAPIROBOOK}) slightly modified, whose proof is straightforward :
\begin{lem}
Let $X\subset Y$ be topological vector spaces, $(\vphi_\lambda)_{\lambda\in\Lambda}$ be a family of continuous operators on $X$ and $Y$. Assume that :
\begin{enumerate}
\item The inclusion is continuous.
\item $X$ is dense in $Y$.
\item $f\in X$ is a common hypercyclic vector for the family $(\vphi_\lambda)_{\lambda\in\Lambda}$, considered as operators on $X$.
\end{enumerate}
Then $f$ is a common hypercyclic vector for the family $(\vphi_\lambda)_{\lambda\in\Lambda}$, considered as operators on $Y$.
\end{lem}
So, if $\alpha>2$, we apply the lemma with :
$$X=L^p\left(\mtr,\frac{dt}{(1+t^2)^2}\right),\quad Y=L^p\left(\mtr,\frac{dt}{(1+t^2)^\alpha}\right).$$
If $1/2<\alpha<2$, set $\veps=\alpha-1/2$. If $f\in\lpmu$, H\"older's inequality gives :
$$\left(\int_\mtr\frac{|f|^p}{(1+t^2)^{1/4+3\varepsilon/4}}\frac{dt}{(1+t^2)^{1/4+\varepsilon/4}}\right)^{1/p}\leq C\left(\int_\mtr\frac{|f|^{2p}}{(1+t^2)^{1/2+3\varepsilon/2}}dt\right)^{1/2p}.$$
Repeated applications of this inequality show that 
$$\exists q\geq 1,\ \exists \beta\geq 2\textrm{ such that }L^q\left(\mtr,\frac{dt}{(1+t^2)^\beta}\right)\subset L^p\left(\mtr,\frac{dt}{(1+t^2)^\alpha}\right),$$
and the lemma works.
\subsubsection{Hyperbolic automorphisms}
We just prove the case $\alpha=2$. First, apply theorem \ref{THMABAKGORDMULT} with the sequence $(\alpha_l)$ identically one, to obtain sequences $(M_k)$ and $(r_k)$. A variant of lemma \ref{LEMSUITEUTILE} gives a nondecreasing sequence $(u_k)$, tending to $+\infty$, and such that :
\begin{enumerate}[a)]
\item $\dis \sum_{k\geq k_0}\frac{u_k}{k^3}<+\infty$.
\item $u_k\sqrt{\frac{r_k}{r_{k-1}}}\xrightarrow{k\to+\infty} 0$.
\item For all $b\in\mtr$, $\dis \sum_{m>k}\frac{u_m}{\left(2^{m-k}\sqrt{\frac{r_k}{r_{k-1}}}-r_kb+b\right)^3}\xrightarrow{k\to+\infty} 0$.
\end{enumerate}
We fix $(f_l)$ a sequence of continuous functions, with $\dis \supp\ f_l\subset\left[-l;-\frac{1}{l}\right]\cup\left[\frac{1}{l};l\right]$, $\dis \left\|f_l\right\|_\infty^p\leq u_l$, and such that, for any $y$ in $\mtr$, for any $\mu\geq 1$,  $\left(f_l(\mu x+\mu y)\right)$ is dense in $\lpdmu$. For $k>k_0$, let us set :
\begin{itemize}
\item $I_k=\left]\frac{1}{\sqrt{r_kr_{k-1}}};\frac{1}{\sqrt{r_kr_{k+1}}}\right[$, and $J_k=-I_k$ : $(I_k)$ and $(J_k)$ are disjoint intervals.
\item $x\in I_k\cup J_k\implies f(x)=f_{j(k)}(r_kx)$. If $x$ is outside $\bigcup_k I_k\cup J_k$, $f(x)=0$.
\end{itemize}
Then $f\in \lpdmu$ : indeed,
\begin{eqnarray*}
\int_0^{+\infty}|f(x)|^pd\mu&\leq&\sum_{k>k_0}\int_{\frac{1}{\sqrt{r_kr_{k-1}}}}^{+\infty}|f_{j(k)}|^pd\mu\\
&\leq&\sum_{k>k_0}u_k\mu\left(\left[\frac{1}{\sqrt{r_kr_{k-1}}};+\infty\right[\right)<+\infty.
\end{eqnarray*}
Fix $\lambda>1$, $b\in\mtr$. We now prove that $f$ is hypercyclic for $S_{\lambda,b}$. Let $l\in\mtn$, $\varepsilon>0$, and $0<\delta<1/2$ whose precise value will be determined later. There exists $k>k_0$ such that :
\begin{itemize}
\item $j(k)=l$.
\item $|\lambda^{M_k}r_k-1|<\delta$.
\item For $m\geq k$, $\dis \sqrt{\frac{r_{m-1}}{r_m}}\geq 2$.
\end{itemize}
Then, 
$$\begin{array}{l}
\dis \int_{-\infty}^{+\infty}|S_{\lambda,b}^{M_k}(f)(x)-f_l(x-b)|^pd\mu\leq\\
\dis \quad\quad\sum_{m<k}\int_{\lambda^{M_k}(x-b)+b\in I_m}\left|f_{j(m)}(\lambda^{M_k}r_m(x-b)+r_mb)-f_l(x-b) \right|^pd\mu\\
\dis \quad\quad+\int_{\lambda^{M_k}(x-b)+b\in I_k}\left|f_{l}(\lambda^{M_k}r_k(x-b)+r_kb)-f_l(x-b) \right|^pd\mu\\
\dis \quad\quad+\sum_{m>k}\int_{\lambda^{M_k}(x-b)+b\in I_m}\left|f_{j(m)}(\lambda^{M_k}r_m(x-b)+r_mb)-f_l(x-b) \right|^pd\mu\\
\dis \quad\quad+S'_1+S'_2+S'_3\\
\leq S_1+S_2+S_3+S'_1+S'_2+S'_3,
\end{array}$$
where $S'_i$ is the same as $S_i$, replacing $I_m$ by $J_m$. Now,
\begin{enumerate}
\item $\dis S_1\leq 2^pu_k\mu\left(\bigcup_{m<k}\frac{I_m-b}{\lambda^{M_k}}+b\right)$. Since
$$\bigcup_{m<k}\frac{I_m-b}{\lambda^{M_k}}+b\subset\left[b-\frac{b}{\lambda^{M_k}};b+\frac{1}{\lambda^{M_k}\sqrt{r_kr_{k-1}}}-\frac{b}{\lambda^{M_k}}\right],$$
we obtain that :
$$S_1\leq 2^{p}u_k\frac{1}{\lambda^{M_k}\sqrt{r_kr_{k-1}}}\leq 2^{p+1}u_k\sqrt{\frac{r_k}{r_{k-1}}}.$$
For $k$ large enough, $|S_1|\leq\varepsilon$.
\item We have 
$\left|\lambda^{M_k}r_k(x-b)+r_kb-(x-b)\right|\leq
\delta|x-b|+r_k|b|$. By uniform continuity of $f_l$, if $\delta$ is small enough, and $k$ is large enough,  $|S_2|\leq\varepsilon$.
\item We have :
\begin{eqnarray*}
S_3&\leq&2^p\sum_{m>k}u_m\mu\left(\frac{I_m-b}{\lambda^{M_k}}+b\right)\\
&\leq&A_1\sum_{m>k}u_m\frac{1}{\left(\frac{r_k}{\sqrt{r_mr_{m-1}}}-r_kb+b\right)^3}\\
&\leq&A_2\dis \sum_{m>k}\frac{u_m}{\left(2^{m-k}\sqrt{\frac{r_k}{r_{k-1}}}-r_kb+b\right)^3},
\end{eqnarray*}
where this last inequality comes from :
$$\frac{r_k}{\sqrt{r_mr_{m-1}}}=\sqrt{\frac{r_k}{r_{k+1}}}\times\dots\times\sqrt{\frac{r_{m-1}}{r_m}}\times\sqrt{\frac{r_k}{r_{m-1}}}\geq 2^{m-k}\sqrt{\frac{r_k}{r_{k-1}}}.$$
For $k$ large enough, $S_3$  is smaller than $\varepsilon$.
\end{enumerate}
$S'_i$ can be treated by the same method as $S_i$ : $f$ is hypercyclic for $S_{\lambda,b}$. Now, as in the holomorphic case, it is not difficult to prove that in fact, for each $\mu\geq 1$, $g(x)=f(\mu x)$ is hypercyclic for all $S_{\lambda,b}$ : this achieves to prove that the set of common hypercyclic vectors is dense!
\section{Final remarks}
\subsection{}Our interest on hypercyclicity  originates from the following \mbox{question :} in \cite{GH}, J.Gordon and H.Hedenmalm characterized the composition operators on the Hilbert space of square summable Dirichlet series $\mch=\{f(s)=\sum_{n\geq 1}a_nn^{-s};\ \|f\|^2:=\sum|a_n|^2<+\infty\}.$ In \cite{BAYHARDIR}, we began a comparison between the properties of the operator $\cphi$ and of its symbol $\phi$. Pursuing this project, we wanted to determinate the hypercyclic composition operators on $\mch$. The answer is very simple :
\begin{thm}
No composition operator on $\mch$ is hypercyclic.
\end{thm}
\begin{pf}
Let $\cphi$ be such a composition operator, induced by $\phi(s)=c_0s+\vphi(s)$, $c_0$ being an integer, and $\vphi$ a Dirichlet series. If $c_0=1$, $\cphi$ is a contraction, and therefore is never hypercyclic. If $c_0=0$, by \cite[lemma 11]{BAYHARDIR}, $\phi_2(\cplus)\subset\mtc_{1/2+\veps}$. Now, take $f$ in $\mch^2$. Then :
$$\left|f\circ\phi_n(+\infty)\right|^2\leq\|f\|^2\zeta(2\Re\phi_n(+\infty))\leq\|f\|^2\max\big(\zeta(1+2\veps),\zeta(2\Re\phi(+\infty))\big).$$
In particular, $\left(\cphi^n(f)\right)$ cannot be dense in $\mch$.
\end{pf}
\subsection{}In \cite{GODSHAPI}, G.Godefroy and H.Shapiro proved that if $\vphi$ is a holomorphic bounded function on $\mtd$, then $M_\vphi^*$ is hypercyclic on $\hdeux$ if and only if $\vphi(\mtd)\cap\mtt\neq\varnothing$. In view of corollary \ref{CORSIMULADJOINTS}, we ask whether there exists a common hypercyclic vector for all $\lambda M_{\vphi^*}$, where $\lambda\vphi(\mtd)\cap\mtt\neq\varnothing$.
\subsection{} In view of theorem \ref{THMSIMULPOIDS}, one may study the weights $\omega$ on $\mtr$ for which the translation operator $Tf(x)=f(x+1)$ and the homothety operator $Sf(x)=f(2x)$ are hypercyclic on $L^1(\mtr,\omega)$.
\begin{defi}
A positive continuois bounded function $\omega$ on $\mtr$ is called a \emph{weight admisible for translation} provided there exists $C>0$ such that, for all $a\in\mtr$, 
$$\int_{a-1}^a\omega(x)dx\leq C\int_{a}^{a+1}\omega(x)dx.$$
It is called \emph{admissible for homothety} if there exists $C>0$ such that, for each $x,y\in\mtr$ with $0\leq x\leq y$ or $x\leq y\leq 0$, 
$$\int_{x/2}^{y/2}\omega(x)dx\leq C\int_x^y\omega(x)dx.$$
\end{defi}
If $\omega$ is admissible for translation (resp. admissible for homothety), the translation operator $T$ (resp. the homothety operator $S$) is continuous on $L^1(\mtr,\omega)$.
\begin{thm}
Let $\omega$ be a continuous bounded positive function on $\mtr$. 
\begin{enumerate}[a)]
\item If $\omega$ is admissible for translation, $T$ is hypercyclic on $L^1(\mtr,\omega)$ if, and only if, there exists $(n_k)_{k\in\mtn}$ a sequence of integers such that :
$$\int_{n_k-q}^{n_k+q}\omega(x)dx\xrightarrow{k\to+\infty}0\textrm{ and }\int_{-n_k-q}^{-n_k+q}\omega(x)dx\xrightarrow{k\to+\infty}0$$
for each $q>0$.
\item If $\omega$ is admissible for homothety, $S$ is hypercyclic on $L^1(\mtr,\omega)$ if, and only if, there exists $(n_k)_{k\in\mtn}$ a sequence of integers such that :
$$\int_{2^{n_k}a}^{2^{n_k}b}\omega(x)dx\xrightarrow{k\to+\infty}0\textrm{ and }\int_{-2^{n_k}b}^{-2^{n_k}a}\omega(x)dx\xrightarrow{k\to+\infty}0$$
for each $0< a\leq b$.
\end{enumerate}
\end{thm}
This statement is the continuous version of Salas theorem \cite{SALASHYPER} on weighted shifts.\\
\begin{pf}
\begin{enumerate}[a)]
\item \textbf{The condition is sufficient :} we apply the hypercyclic criterion. Let $(P_j)$ be a dense sequence in $L^1(\mtr,\omega)$ of bounded functions compactly supported. If $\supp P_j\subset[-q,q]$, then :
\begin{eqnarray*}
\left\|T^{n_k}P_j\right\|&=&\int_{-n_k-q}^{-n_k+q}\left|P_j(x+n_k)\right|\omega(x)dx\\
&\leq&\|P_j\|_{\infty}\int_{-n_k-q}^{-n_k+q}\omega(x)dx\xrightarrow{k\to+\infty}0.
\end{eqnarray*}
Take $Af(x)=f(x-1)$. $A$ is a (possibly unbounded) right inverse of $T$. It is straightforward that $\|A^{n_k}P_j\|\xrightarrow{k\to+\infty}0$.\\
\textbf{The condition is necessary :} By a diagonal argument, it suffices to prove that, for all $\veps>0$, for all $q>0$, there exists $N$ arbitrarily large such that 
$$\int_{N-q}^{N+q}\omega(x)dx\leq\veps\textrm{ and }\int_{-N-q}^{-N+q}\omega(x)dx\leq\veps.$$
We set $A_1=\inf_{[-q,q]}\omega$, $A_2=\sup_\mtr\omega$. Since the set of hypercyclic vectors for $T$ is dense, there is a hypercyclic vector $f\in L^1(\mtr,\omega)$ such that :
\begin{eqnarray}\label{EQWEIGHT1}
\left\|f-1_{[-q,q]}\right\|\leq\frac{\veps A_1}{2A_2}.
\end{eqnarray}
We can also find $N$ arbitrarily large, $N>2q$, such that :
\begin{eqnarray}\label{EQWEIGHT2}
\left\|T^Nf-1_{[-q,q]}\right\|\leq\frac{\veps A_1}{2A_2}.
\end{eqnarray}
Since $N\geq 2q$, inequality (\ref{EQWEIGHT1}) implies 
$$\int_{N-q}^{N+q}|f(x)|\omega(x)dx\leq\frac{\veps}{2},$$
whereas inequality (\ref{EQWEIGHT2}) gives :
$$\int_{-q}^q|f(x+N)-1|\omega(x)dx\leq\frac{\veps A_1}{2A_2},$$
which in turn proves :
$$\int_{N-q}^{N+q}|f(x)-1|\omega(x)dx\leq\frac{\veps}{2}.$$
Thus
$$\int_{N-q}^{N+q}\omega(x)dx\leq\veps.$$
We proceed with the same method for the other inequality.
\item We prove that the condition is sufficient by using another time Kitai's criterion, now with continuous functions whose supports are contained in intervals like $[-A,-\delta]\cup[\delta,A]$, with $0<\delta\leq A$. For the necessity, we fix $0<a\leq b$, and $A_1=\inf_{[a,b]}\omega$, $A_2=\sup_\mtr\omega$. There exists $f\in L^1(\mtr,\omega)$ and $N$ arbitrarily large (in particular, $2^Na>b$) with :
\begin{eqnarray}\label{EQWEIGHT3}
\left\|f-1_{[a,b]}\right\|\leq\frac{\veps A_1}{2A_2},
\end{eqnarray}
\begin{eqnarray}\label{EQWEIGHT4}
\left\|S^Nf-1_{[a,b]}\right\|\leq\frac{\veps A_1}{2A_2}.
\end{eqnarray}
As previously, (\ref{EQWEIGHT3}) gives :
$$\int_{2^Na}^{2^Nb}|f(x)|\omega(x)dx\leq\frac{\veps}{2},$$
and (\ref{EQWEIGHT4}) implies :
$$\int_{2^Na}^{2^Nb}|f(x)-1|\omega(x)dx\leq\frac{\veps}{2}.$$
This in turn implies 
$$\int_{2^Na}^{2^Nb}\omega(x)dx\leq\veps.$$
\end{enumerate}
\end{pf}
\begin{ex}
For the weight $\dis\omega(x)=\frac{1}{1+|x|}$, the translation operator $T$ is hypercyclic, whereas the homothety operator $S$ is not.
\end{ex}

\end{document}